# Analytical solution of transient scalar wave and diffusion problems of arbitrary dimensionality and geometry by RBF wavelet series

W. Chen[*]

Department of Informatics, University of Oslo, P.O.Box 1080, Blindern, 0316 Oslo, Norway

Email: wenc@ifi.uio.no

**Summary**

This study applies the RBF wavelet series to the evaluation of analytical solutions of linear time-dependent wave and diffusion problems of any dimensionality and geometry. To the best of the author's knowledge, such analytical solutions have never been achieved before. The RBF wavelets can be understood an alternative for multidimensional problems to the standard Fourier series via fundamental and general solutions of partial differential equation. The present RBF wavelets are infinitely differential, compactly supported, orthogonal over different scales and very simple. The rigorous mathematical proof of completeness and convergence is still missing in this study. The present work may open a new window to numerical solution and theoretical analysis of many other high-dimensional time-dependent PDE problems under arbitrary geometry.

---

[*] This research is supported by Norwegian Research Council.

# 1. Fourier series, a historic retrospect

Many of the important concepts of analysis and computation have their origins in the study of physical problems leading to the partial differential equation (PDE) systems [1]. The currently ubiquitous Fourier series and transform came from Fourier's original exploration of the solution of a bar heat transmission problem in the early 1800s. What Fourier proposed due to this quest is that an arbitrary 1D function $f(x)$ over a bounded interval, even if not differentiable, can be represented by an infinite sum of sinusoids

$$f(x) = \sum_{k=1}^{\infty} c_k \sin \frac{k\pi}{a} x, \qquad (1)$$

where sinusoids are the eigenfunctions of this PDE problem. Despite a lack of rigorous proof, Fourier was quite confident of the basic truth of his assertion for obvious physical and geometric grounds. Nevertheless, the implications of this discovery go well beyond Fourier's wildest imaginations. The unanswered mathematical points forcefully gave birth to many more new problems and consequently motivated the development of many important mathematic concepts and techniques such as Riemann integration, Sturm-Liouville eigenvalue problem, set theory, Laplace transform, Lebesgue integration, Green's function and distribution theory, functional analysis, and most recently, wavelets theory [1] as well as enormous applications in numerous ramifications of science and engineering.

Despite the widespread applicability [1-3], the Fourier analysis approach suffers some drawbacks. Most noticeably, for more than one-dimension problems, the direct use of the Fourier series becomes very mathematically complicated and is only feasible for such regular geometry as rectangular, circle, sphere, cylindrical domains, etc., where we can separate the space variables (in Cartesian, polar, or some other coordinate systems [4]). Otherwise the tensor product approach, very costly in high dimensions, must be applied. However, when the scattered data are involved, the tensor product Fourier analysis also immediately fails.

Because of the great success of the polynomials, splines, and tensor product methods, mathematician and engineers alike grow accustomed to expressing a function in terms of coordinate variables. To majority of scientific and engineering community, the radial basis function (RBF), which uses the one-dimensional distance variable irrespective of dimensionality, is a quite brand-new and exotic concept. In high dimensional scattered data cases, the RBF approach, however, is the method of the choice [5]. It is also found that the RBF is very efficient in handling lower-dimensional problems [5-7]. In the following sections, we will try to clarify an underlying connection between the RBF and PDE problem, which unveils that the RBF approach, just like Fourier series, is no more than a natural technique originating from the solution of PDE problems.

Another notorious drawback of the Fourier approach is to lack efficient localization in both time and frequency or scale, where the promise of wavelets technique comes into play. In recent decade fast development and widespread applications of the wavelets have been one of the most significant achievements in mathematics and many physical areas [8]. Unlike the most traditional expansion systems such as Fourier one [9], the basis functions of the wavelet analysis, however, are not solutions of differential equation. This comes without a surprise since the wavelets have mainly been developed in signal processing area, where the PDE system is rarely involved. By means of the fundamental solution and general solution of partial differential equation, the next section concerns a convergence of the RBF and wavelets while attainting the compactness and infinite smoothness.

## 2. Radial basis function and wavelets

Since pioneer works of Frankle [10], Michaelli [11] and Kansa [7], the research into the RBF theory and applications becomes very active. In parallel, Daubechies' breakthrough orthogonal compact wavelets [12] lead revolutionary advances in multiscale analysis. The RBF is well known for its striking effectiveness in multivariate scattered data approximation [5]. However, in general the RBFs available now lack critical multiscale analysis capability. To handle high-dimensional multiscale scattered data and PDE systems, the RBF wavelets are mostly wanted to combine the strengths of both. In last decade much effort has been devoted to various non-orthogonal prewavelets RBF theory and applications by using some constructive approximation strategies [13-17]. Very recently Chen [18] could develop the orthogonal RBF wavelet series and transform by using the fundamental solution and general solution of some typical PDE's. The work given in [18], however, is more conjecture and speculation than a complete theory with some evident errors. Next is summarized the discrete orthogonal RBF wavelet series versus the Fourier series. Some corrections to ref. 18 are also given.

Consider a real-value Riemann-integral multidimensional function $f(\bar{x})$ on the $n$-dimension spherical domain $\Omega$ of radius $R$,

$$f(\bar{x}) = a_0/2 + \sum_{j=1}^{\infty} \sum_{k=1}^{\infty} \alpha_{jk} \varphi_n\left(\eta_j \|\bar{x} - \bar{x}_k\|\right) \tag{2}$$

is called its RBF-based wavelet expansion series, where $\varphi_n$ represents the wavelet basis function, $\|\bar{x} - \bar{x}_k\|$ means Euclidean distance; and $\alpha_{jk}$ are the expansion coefficients. The reason why $\alpha_0$ is multiplied by ½ is technical as in the Fourier series and will be explained below. To discrete harmonic analysis, we choose the nonsingular general solutions of *n*-dimension Helmholtz equation

$$\varphi_n(r_k) = 1, \quad \eta = 0, \tag{3a}$$

$$\varphi_1(\eta_j r_k) = \frac{1}{2\eta_j}\sin(\eta_j r_k), \quad \eta \neq 0, \tag{3b}$$

$$\varphi_n(\eta_j r_k) = \frac{1}{4}\left(\frac{\eta_j}{2\pi r_k}\right)^{(n/2)-1} J_{(n/2)-1}(\eta_j r_k), \quad n \geq 2, \quad \eta \neq 0 \tag{3c}$$

as wavelet basis function, where $r_k = \|\bar{x} - \bar{x}_k\|$, $J_{(n/2)-1}()$ is the (*n*/2-1)th order Bessel function of the first kind; and $\eta_j$ are the zeros of $\varphi_n(Rr)$. The eigenfunctions (3) form an orthonormal set of basis functions over different scales. In 3D case, the present RBF is the renowned Sinc function.

Despite mathematical beauty and simplicity of the RBF expansion (2), a theoretical proof of convergence and completeness under certain conditions is missing in [18]. Levesley et al [19] provides an important clue to integrate the RBF wavelet approach with convolution operator theory. Considering the fundamental solution in the distribution theory, the general solution of *n*-dimensional Helmholtz equation may (may not) establish

$$\int_{\Omega_\infty} \varphi_n(\|\bar{x} - \bar{x}_k\|) d\Omega_k \neq 0 \prec \infty, \tag{4}$$

and then for a radially symmetric domain

$$\int_0^\infty r^{n-1} \varphi_n(r) dr \prec \infty. \tag{5}$$

Note that the above condition is very stringent. Thus, for $f(\bar{x}) \in R^n$ we have

$$\lim_{j \to \infty} c_j \int_\Omega \varphi_n(j\|\bar{x} - \bar{x}_k\|) f(\bar{x}) d\Omega_k = f(x), \tag{6}$$

where $c_j$ are the coefficients depending on the RBF $\varphi_n$. After a numerical discretization of the integral in Eq. (6), we get RBF approximate expression (2). Its convergence is therefore guaranteed uniform and compact [19]. The condition (4) is satisfied by RBF wavelet transform using singular fundamental solution given in later section 6. Like the Fourier expansion, the function expressible in RBF wavelet series (2) should subject to some conditions, especially if $\varphi_n$ does not satisfy the condition (4). For now, we do not have an explicit answer for this. A more loose condition on convergence of the Helmholtz RBF wavelet series may exist.

On the other hand, very recently the present author became aware of a relevant important paper by Mourou and Trimeche [20], which originally introduces the so-called generalized wavelet transform using the solution of the Bessel operator. Although Mourou and Trimeche [20] use the term of the radial function, the generalized wavelet transform given there is in fact based on 1D coordinate variable rather than Euclidean distance variable and therefore is in fact not a RBF wavelet, which can only handle one-dimensional cases if without tensor product approach. Noteworthily, the admissibility condition of their generalized wavelet transform is also applicable to the present RBF wavelet series due to the underlying connection between Helmholtz equation and Bessel operator. Namely, we have

$$0 \prec C_g = \int_0^\infty |\hbar(\varphi_n)(\lambda)|^2 \frac{d\lambda}{\lambda} \prec \infty, \tag{7}$$

where $\hbar(\varphi_n)$ is the Fourier transform of general solution or fundamental solution $\varphi_n$ of Helmholtz equation.

Another important issue is how to calculate the expansion coefficients. It is known that the Bessel function possesses the orthogonality:

$$\int_0^R r J_v(\eta r) J_v(\mu r) dr = \begin{cases} 0 & \eta \neq \mu \\ R^2 J_{v+1}(R\eta)^2 / 2 & \eta = \mu \end{cases} \quad (8)$$

where $J_v(R\eta)=J_v(R\mu)=0$. However, the present RBF wavelet series is not necessary orthogonal over translation of the same scale. Ref. [18] made wrong assumption that the present RBF wavelets are orthogonal over both scale and translation. In this regard the Gram-Schmidt orthogonality method may enforce the present RBF wavelets orthogonal over translation. Thus we can directly determine the expansion coefficients.

In the 1D case, the coefficient formulas are quite similar to those of the corresponding Fourier series:

$$\alpha_0 = \frac{2}{R} \int_{\Omega_R} f(\zeta) d\Omega_\zeta, \quad (9a)$$

$$\alpha_{jk} = \frac{2j\pi}{R^2} \int_{-R}^{R} f(\zeta) \sin\left(\frac{j\pi}{R} r_{k\zeta}\right) d\zeta, \quad (9b)$$

$$j=1,2,\ldots, \quad k=1,2,\ldots$$

For multidimensional problems, we have

$$\alpha_0 = \frac{2}{V_n} \int_{\Omega_R} f(\zeta) d\Omega_\zeta, \quad (10a)$$

$$\alpha_{jk} = \frac{8}{S_n J_{n/2}(\eta_j)^2} \left(\frac{\eta_j}{2\pi}\right)^{1-n/2} \int_{\Omega_\zeta} r_{k\zeta}^{1-n/2} f(\zeta) J_{(n/2)-1}(\eta_j r_{k\zeta}) d\Omega_\zeta, \tag{10b}$$

$$j = 1,2,\ldots, \quad k = 1,2,\ldots, \quad n \geq 2,$$

where $r_{k\zeta} = \|\bar{x}_k - \bar{x}_\zeta\|$, $R$ is the radius of the spherical domain centering node $k$, and $S_n$ represents the surface parameter; $V_n$ is the volume of $n$-dimensional sphere. The reason for the coefficient ½ in (2) is that this will facilitate to streamline formulas (9a) and (9b) or (10a) and (10b) if we choose proper constant in Eq. (3a).

### 3. Analytical solution of transient wave problem with RBF wavelets

The following $n$-dimension homogeneous wave problem with homogeneous boundary conditions serves as an illustrative example of the present strategy:

$$\nabla^2 u = \frac{1}{c^2} \frac{\partial^2 u}{\partial t^2}, \qquad x \in \Omega, \tag{11}$$

$$\begin{cases} u(\bar{x},t) = 0, & \bar{x} \subset S_u, \\ \dfrac{\partial u(\bar{x},t)}{\partial n} = 0, & \bar{x} \subset S_T, \end{cases} \quad t \geq 0, \tag{12}$$

$$\begin{cases} u(\bar{x},0) = \phi(\bar{x}), & \bar{x} \subset S, \\ u_t(\bar{x},0) = \psi(x) & \bar{x} \subset S, \end{cases} \tag{13}$$

where $\bar{x}$ means multi-dimensional independent variable, and $n$ is the unit outward normal. $S = S_u \cup S_\Gamma$, The separation of the time variable yields

$$u(\bar{x},t) = T(t) v(\bar{x}) \tag{14}$$

Substitution into the wave equation (11) allows the separation of the time part

$$\nabla^2 v + \lambda^2 v = 0, \tag{15}$$

$$\frac{d^2 T}{dt^2} + \lambda^2 T = 0. \tag{16}$$

Here $\lambda$ is the separation constant and actually eigenvalue of the system. The problem has only nonnegative eigenvalues $\lambda_j$ [22, pp. 248] and the solutions of the wave equation are

$$u(\bar{x},t) = \frac{1}{2} A_0 + \frac{1}{2} B_0 t + \sum_{j=1}^{\infty} \left[ A_j \cos\left(ct\sqrt{\lambda_j}\right) + B_j \sin\left(ct\sqrt{\lambda_j}\right) \right] v_j(\bar{x}), \tag{17}$$

where $v_j(\bar{x})$ is the corresponding eigenfunctions. For the Robin (radiation) boundary condition ($\partial u/\partial n + au = 0$), the nonnegative eigenvalues holds provided that $a \geq 0$ [22].

**3.1. Helmholtz eigenvalues and eigenfunctions with RBF**

The challenging issues are to evaluate the spatial eigenvalues and eigenfunctions from the Helmholtz equation (15). The standard approaches based on the coordinate variables do not work in each case of irregular geometry, scattered data and high dimensions. So we have to resort the radial basis function methodology. In the case of the *n*-dimensional Helmholtz problem with eigenvalues $\lambda_j$, the complete RBF solution contains a complex argument [23]:

$$u^*(r) = \frac{1}{4}\left(\frac{\lambda}{2\pi r}\right)^{n/2-1} \left[ Y_{n/2-1}(\lambda r) - iJ_{n/2-1}(\lambda r) \right], \quad n \geq 2, \tag{18}$$

where $r_k = \|x - x_k\|$, $J_{(n/2)-1}$ and $Y_{(n/2)-1}$ are respectively the ($n/2$-1)th order Bessel function of the first and second kinds. The Bessel function of the first kind is $C^\infty$ smooth while the Bessel function of the second kind encounters a singularity at the origin. Thus, the nonsingular general solution (regular distribution) is the imaginary part of the above complex RBF solution (18) and is seen as the eigenfunctions. Namely, we have RBF eigenfunctions as follows

$$\varphi_n(r_k) = 1, \quad \eta = 0, \tag{19a}$$

$$\varphi_1(\lambda_j r_k) = \frac{1}{2\lambda_j} \sin(\lambda_j r_k), \quad n=1, \ \lambda \neq 0, \tag{19b}$$

$$\varphi_n(\lambda_j r_k) = \frac{1}{4}\left(\frac{\lambda_j}{2\pi r_k}\right)^{(n/2)-1} J_{(n/2)-1}(\lambda_j r_k), \quad n \geq 2, \ \lambda \neq 0, \tag{19c}$$

Since the eigenfunctions must be finite at the origin, scraping the singular part of the RBF solution (18) does not raise the completeness issue in general.

Note that the general solution is independent of geometry. In other words, the general solution satisfies the Helmholtz equation (15) irrespective of the boundary shape of interests, which is validated by the computer software "Maple" in terms of 2D and 3D Cartesian coordinates.

The next issue is to calculate the eigenvalues. We present the two boundary knot method (BKM) [24,25] schemes for this task. The first scheme is to employ the BKM discretization of the Helmholtz equation (15). The approximate representation of the symmetric BKM is

$$v(x) = \sum_{s=1}^{L_D} \beta_s \varphi_n(r_s) - \sum_{s=L_D+1}^{L_D+L_N} \beta_s \frac{\partial \varphi_n(r_s)}{\partial n}, \tag{20}$$

where $k$ is the index of source points on boundary, $\beta_k$ are the desired coefficients; $n$ is the unit outward normal as in boundary condition (12), and $L_D$ and $L_N$ are respectively the numbers of knots on the Dirichlet and Neumann boundary surfaces. The minus sign associated with the second term is due to the fact that the Neumann condition of the first order derivative is not self-adjoint.

In terms of representation (20), we have the homogeneous collocation analogue $Hv=0$ of boundary condition equations (12), where $H$ is the symmetric BKM interpolation matrix. And then, just as in the traditional Fourier solution of 1D wave problem, the determinant of interpolation matrix has to be zero to attain the nontrivial solution, namely,

$$\det(H)=0. \qquad (21)$$

The infinitely many roots of the above algebraic equation are the eigenvalues of the Helmholtz equation (15). Note that since we use the symmetric BKM [25], all the solutions of eigenvalues will be real valued.

The solution of transcendental equation (21) is often a daunting job and hence less attractive in practical use. So we develop the second strategy of the BKM for this task. The Helmholtz equation (15) can be rewritten as

$$\nabla^2 v + \delta^2 v = \left(-\lambda^2 + \delta^2\right)v, \qquad (22)$$

where $\delta$ is a small artificial real parameter (around 0.1) and insensitive to the boundary shape and dimensionality [25]. In terms of the symmetric BKM expression (20), the above equation (22) can be discretized into the standard algebraic eigenvalue problem

$$Kv = \left(-\lambda^2 + \delta^2\right)v, \qquad (23)$$

where the interpolation matrix $K$ is symmetric irrespective of the boundary conditions. Many algorithm packages are readily available to calculate this standard eigenvalue problem. Unlike the previous first strategy, the second scheme, however, requires using some inner nodes to guarantee the stability and accuracy of the BKM solution [25].

**3.2. Analytical solution with RBF wavelet series**

In terms of solution (17), the resulting analytical solution of wave equations (1,2,3) can be expressed as

$$u(\bar{x},t) = \frac{1}{2}A_0 + \frac{1}{2}B_0 t + \sum_{j=1}^{\infty}\sum_{k=1}^{\infty}\left[A_{jk}\cos\left(ct\sqrt{\lambda_j}\right) + B_{jk}\sin\left(ct\sqrt{\lambda_j}\right)\right]\varphi_n\left(\lambda_j\|\bar{x}-\bar{x}_k\|\right). \quad (24)$$

The RBF series solution (24) is valid for any dimensionality and geometry since the RBF approach is independent of dimensionality and geometry. Note that this solution is an inseparable wavelet series, where eigenvalue $\lambda_j$ is understood scale parameter (dilation) and source node $\bar{x}_k$ are seen as location parameter (shift) in wavelet terminology.

In terms of initial conditions (13), we have

$$u(\bar{x},0) = \frac{1}{2}A_0 + \sum_{j=1}^{\infty}\sum_{k=1}^{\infty} A_{jk}\varphi_n\left(\lambda_j\|\bar{x}-\bar{x}_k\|\right) = \phi(\bar{x}), \quad (25a)$$

$$u_t(\bar{x},0) = \frac{1}{2}B_0 + \sum_{j=1}^{\infty} B_{jk}c\sqrt{\lambda_j}\sum_{k=1}^{\infty}\varphi_n\left(\lambda_j\|\bar{x}-\bar{x}_k\|\right) = \psi(\bar{x}). \quad (25b)$$

$\phi(\bar{x})$ and $\psi(\bar{x})$ are expressible in the RBF wavelet series (2) of Helmholtz general solution provided that in this case they have the first order differential continuity and equal zero at boundary. Like the Fourier series approach, the coefficients $A_j$ and $B_j$ can most efficiently calculated via orthogonality. Applying the divergence theorem and Green's second identity, it is proved [22, pp. 246] that for Helmholtz equation (15), the real eigenfunctions corresponding to distinct real eigenvalues are necessarily orthogonal. Concerning the eigenfunctions (19), the orthogonality over scale is

$$\int_\Omega X_{sp} X_{tq} d\Omega = 0 , \quad \text{for } s \neq t , \tag{26}$$

where *s* and *t* denote scales, *p* and *q* are the translation locations. It is noted that the RBF wavelet series is not necessarily orthogonal over the translation. By the Gram-Schmidt orthogonality method, we can enforce orthogonization of general solution eigenfunctions of the same scale. Otherwise, we have to solve the simultaneous equations for coefficients $A_j$ and $B_j$.

In [22, pp. 254-256] the method of separation of variables is used to solve vibrations of a regular drumhead, where the eigenfunctions using coordinate variable is a general Fourier series combining Bessel function and trigonometric functions. Due to the orthogonality, all expansion coefficients, the number of which is equal to that of the present RBF wavelet coefficients, are calculated directly. By analogy with this strategy, the following formulas are given without rigorous proof under translation-orthogonal assumption and just for reference. A relevant Gram-Schmidt procedure should be developed in the future.

$$A_0 = \frac{2}{V_n} \int_{\Omega_R} \phi(\zeta) d\Omega_\zeta , \tag{27}$$

$$B_0 = \frac{2}{V_n c \sqrt{\lambda_j}} \int_{\Omega_R} \psi(\zeta) d\Omega_\zeta , \tag{28}$$

$$A_{jk} = \frac{2}{S_n J_{n/2}(\lambda_j)^2} \left(\frac{\lambda_j}{2\pi}\right)^{1-n/2} \int_{\Omega_\zeta} r_{k\zeta}^{1-n/2} \phi(\zeta) J_{(n/2)-1}(\lambda_j r_{k\zeta}) d\Omega_\zeta , \tag{29}$$

$$B_{jk} = \frac{8}{c \sqrt{\lambda_j} S_n J_{n/2}(\lambda_j)^2} \left(\frac{\lambda_j}{2\pi}\right)^{1-n/2} \int_{\Omega_\zeta} r_{k\zeta}^{1-n/2} \psi(\zeta) J_{(n/2)-1}(\lambda_j r_{k\zeta}) d\Omega_\zeta , \tag{30}$$

$$j = 1,2,\ldots, \quad k = 1,2,\ldots$$

It is very interesting to note that the 1D nonsingular general solution (19b) is the same sinusoids as in the 1D Fourier series solution (1). The distinctness of the RBF wavelet and Fourier series solutions is that when applied to practical problems with a truncated finite series, the former can locally adjust the scale parameter to avoid the Gibbs phenomena. Therefore, the RBF wavelet series solution are much more robust than the Fourier one.

It is worth pointing out that the present RBF wavelets using the general solution of Helmholtz equation have periodic (harmonic) property and are much natural than those periodic RBF developed by [26,27], where the sine and cosine functions are employed as the RBF basis without considering dimensional affect and PDE eigensolution.

## 4. Applications to inhomogeneous problems

In practical engineering, the governing equation [11] and boundary conditions [12] are often not homogeneous. For instance,

$$\nabla^2 u = \frac{1}{c^2}\frac{\partial^2 u}{\partial t^2} + f(\bar{x}), \qquad x \in \Omega \qquad (31)$$

$$\begin{cases} u(\bar{x}) = D(\bar{x}), & \bar{x} \subset S_u, \\ \frac{\partial u(\bar{x})}{\partial n} = R(\bar{x}), & \bar{x} \subset S_T, \end{cases} \qquad t \geq 0, \qquad (32)$$

$$\begin{cases} u(\bar{x},0) = \phi(\bar{x}), & \bar{x} \subset S, \\ u_t(\bar{x},0) = \psi(x) & \bar{x} \subset S, \end{cases} \qquad (33)$$

where $f(\bar{x})$ is outside forcing function. Just like the Fourier solution of inhomogeneous PDE's, we have two strategies to apply the method of the separation of variables [22, pp.

140-143] with RBF wavelet series. The method of shifting the data makes the inhomogeneous boundary conditions homogeneous by subtracting any known function that satisfies them, while the expansion method expand everything in the eigenfunctions of the corresponding homogeneous problems and then we can get a set of ordinary differential equations which is very easy to solve.

## 5. Applications to other problems

Besides the wave equation considered previously, the following equations also frequently appear in electrical, magnetic, thermal, gravitational, vibration, hydrodynamics and acoustics problems [4]:

1. The diffusion equation

$$\nabla^2 u = \frac{1}{h^2}\frac{\partial u}{\partial t}, \tag{34}$$

2. The damped wave equation

$$\nabla^2 u = \frac{1}{c^2}\frac{\partial^2 u}{\partial t^2} + R\frac{\partial u}{\partial t}, \tag{35}$$

3. Transmission line equation

$$\nabla^2 u = \frac{1}{c^2}\frac{\partial^2 u}{\partial t^2} + R\frac{\partial u}{\partial t} + Su, \tag{36}$$

The solution of all these equations can be reduced to solutions of the scalar Helmholtz equation (15) and the corresponding very simple ordinary differential equation in time. Thus, the extension of the present RBF wavelet series solution to these equations is very straightforward and omitted here for brevity.

## 6. Generalized RBF wavelet series and transforms

In the foregoing sections, we only consider the RBF wavelet series using the general solution of Helmholtz equation. Those RBF wavelets describe periodic behaviors of physical systems as the trigonometric Fourier series. The higher order fundamental solutions of Laplace operator can construct the RBF wavelet series corresponding to the common polynomial interpolation and approximation, where the order of fundamental solution corresponds to the scale or polynomial order. Moreover, they are orthogonal over scales since the lower-order fundamental solutions satisfy the higher-order operator. The drawbacks of the Laplacian RBF wavelets are limited smoothness.

On the other hand, we can create continuous RBF wavelet transform using the fundamental solution (irregular distribution) of PDE's. For example,

$$F(\lambda,\xi) = \int_{\Omega_\infty} f(\zeta)\overline{g_n(\lambda r_{\xi\zeta})} d\Omega_\zeta, \tag{37}$$

and

$$f(\vec{x}) = C_g^{-1} \int_0^{+\infty} \int_{\Omega_\infty} F(\lambda,\xi) g_n(\lambda r_{x\xi}) \lambda^{2n-1} d\Omega_\xi d\lambda, \tag{38}$$

where $C_g$ is decided by formula (7),

$$g_n(\lambda r_k) = \frac{1}{2\pi}\left(\frac{-i\lambda}{2\pi r_k}\right)^{(n/2)-1} K_{(n/2)-1}(-i\lambda r_k) \tag{39}$$

for Helmholtz harmonic wavelets and

$$g_n(\lambda r_k) = \frac{1}{2\pi}\left(\frac{\lambda}{2\pi r_k}\right)^{(n/2)-1} K_{(n/2)-1}(\lambda r_k) \tag{40}$$

for modified Helmholtz wavelets, $K$ is the modified Bessel function of the first kind. The corresponding dual wavelet basis function is its conjugate function $\overline{g_n(\lambda r_k)}$, which satisfies the condition (4). The modified Helmholtz wavelets are the counterpart of the Laplace transform. Note that here we correct some errors of Eqs. (15) and (16) in ref. 18 without using $r^{n-1}$ weight. For more details see refs. [18,21]. The harmonic RBF wavelet transform enjoys a nice feature of Fourier transform:

$$\nabla^2 F(\lambda,\xi) = -\lambda^2 F(\lambda,\xi). \tag{41}$$

It is also highly likely to develop the RBF wavelets transforms based on the general and fundamental solutions of the other typical partial differential equation such as convection-diffusion equation [18].

### 7. Promises and open problems

**Promises:**

1. Since the present series solution is wavelets, the Gibbs phenomenon long bothering the Fourier series is eliminated. By adapting the arbitrarily scaling parameter (dilations) rather than the dyadic multiresolution analysis, we get locally supported RBF basis both in scale and location. Compactly-supported wavelets using spline, especially orthogonal such as the popular Daubechies wavelets [12], have a limited degree of smoothness in compromise to the compactness. In contrast, the RBF wavelets are not only orthogonal over scale but also infinitely differentiable.
2. Due to the use of the RBF, we avoid using the tensor-product approach for high dimensional problems with irregular geometry. In addition, the present method is meshfree and feasible to handle scattered data problems.
3. The strategy presented here is expected feasible to other type of problems such as data processing and edge detection etc., where an efficient description of multidimensional multiscale scattered data is crucial. In addition, the present method may be employed to evaluate the fractional differential equation corresponding to fractal geometry by using fractional dimensionality.

**Problems**:

The present form of the RBF wavelets, resembling the immature status of Fourier's early work, which used concepts and theories as yet undeveloped or underdeveloped [1], may provide some of the refresh impulse on new advances in applied and basic analysis. Some worrisome points of this study are stated below (the distribution theory may be useful to research them).

1. In terms of Fourier series, there are three kinds of convergence: pointwise, uniform and $L^2$ convergences. This study will not go further into the convergence issue of the RBF wavelet series solution of transient PDE's. Ref. [22] defines that an orthogonal system is called complete if and only if it is not a proper subset of another orthogonal system. The present general solution eigenfunctions satisfy this condition. In addition, the key issue is if it satisfies Parseval's identity. Section 2 provides some proof of the convergence and completeness of the RBF wavelet series via convolution operator theory and admissibility condition of wavelets. The proof of completeness and convergence of RBF series solution of PDE's, however, is still lacking now. Chen and Tanaka [17] and Hon and Chen [28] also physically discussed the completeness issue of the general solution expansion series within the framework of the boundary knot method.
2. It is unknown now if the RBF wavelet transform given in section 6 [18] can be employed to get the analytical solution of time-dependent problems in the integral form. If OK, this integral solution should be equal to the limit of the present RBF wavelet series solution.
3. The computing formulas (9), (10) and (27-30) of expansion coefficients based on orthogonality may be error prone.